\documentclass[a4paper,12pt]{amsart}
\usepackage{amsfonts}
\usepackage{amssymb}
\usepackage{ifthen}
\usepackage{graphicx}
\usepackage{enumitem}
\usepackage{amsmath}
\nonstopmode \numberwithin{equation}{section}
\usepackage[usenames]{color}
\newtheorem{thm}{Theorem}[section]
\newtheorem{lem}{Lemma}[section]
\newtheorem{cor}{Corollary}[section]
\newtheorem{prop}[thm]{Proposition}

\theoremstyle{definition}
\newtheorem{mlem}{Main lemma}[section]
\newtheorem{assertion}{Assertion}[section]
\newtheorem{cl}{Claim}[section]
\newtheorem{ca}{Case}[section]
\newtheorem{sca}{Subcase}[section]
\newtheorem{scl}{Subclaim}[section]
\newtheorem{conj}[thm]{Conjecture}
\newtheorem{fact}{Fact}[section]
\newtheorem{defn}{Definition}[section]
\newtheorem{prob}{Problem}[section]
\newtheorem{ques}[thm]{Question}
\newtheorem{rem}{Remark}[section]
\newtheorem{exam}{Example}[section]

\numberwithin{equation}{section}

\newcounter {own}
\def\theown {\thesection       .\arabic{own}}

\newlist{steps}{enumerate}{1}
\setlist[steps,1]{
  leftmargin=*,
  label=\textbf{Step \arabic*}.,
  ref=Step~\arabic*,
}

\newenvironment{pf}[1][]{%
 \vskip 3mm
 \noindent
 \ifthenelse{\equal{#1}{}}%
  {{\slshape Proof. }}%
  {{\slshape #1.} }%
 }%
{\qed\bigskip}

\newcounter{alphabet}
\newcounter{tmp}
\newenvironment{Thm}[1][]{\refstepcounter{alphabet}%
\bigskip%
\noindent%
{\bf Theorem \Alph{alphabet}}%
\ifthenelse{\equal{#1}{}}{}{ (#1)}%
{\bf .} \itshape}{\vskip 8pt}

\newcounter{alphabet2}

\makeatletter
\newcommand{\Ref}[1]{\@ifundefined{r@#1}{}{\setcounter{tmp}{\ref{#1}}\Alph{tmp}}}
\makeatother

\newenvironment{Lem}[1][]{\refstepcounter{alphabet}%
\bigskip%
\noindent%
{\bf Lemma \Alph{alphabet}}%
{\bf .} \itshape}{\vskip 8pt}

\DeclareMathOperator*{\esssup}{ess\,sup}

\newcommand{\ID}{{\mathbb D}}

\def\be{\begin{equation}}
\def\ee{\end{equation}}

\newcommand{\ben}{\begin{enumerate}}
\newcommand{\een}{\end{enumerate}}

\newcommand{\blem}{\begin{lem}}
\newcommand{\elem}{\end{lem}}
\newcommand{\bthm}{\begin{thm}}
\newcommand{\ethm}{\end{thm}}
\newcommand{\bcor}{\begin{cor}}
\newcommand{\ecor}{\end{cor}}
\newcommand{\beg}{\begin{exam}}
\newcommand{\eeg}{\end{exam}}
\newcommand{\begs}{\begin{examples}}
\newcommand{\eegs}{\end{examples}}
\newcommand{\bdefe}{\begin{defn}}
\newcommand{\edefe}{\end{defn}}
\newcommand{\bprob}{\begin{prob}}
\newcommand{\eprob}{\end{prob}}
\newcommand{\bques}{\begin{ques}}
\newcommand{\eques}{\end{ques}}
\newcommand{\bei}{\begin{itemize}}
\newcommand{\eei}{\end{itemize}}
\newcommand{\bcon}{\begin{conj}}
\newcommand{\econ}{\end{conj}}
\newcommand{\bop}{\begin{op}}
\newcommand{\eop}{\end{op}}

\newcommand{\bas}{\begin{assertion}}
\newcommand{\eas}{\end{assertion}}

\newcommand{\bfa}{\begin{fact}}
\newcommand{\efa}{\end{fact}}

\newcommand{\bca}{\begin{ca}}
\newcommand{\eca}{\end{ca}}
\newcommand{\bsca}{\begin{sca}}
\newcommand{\esca}{\end{sca}}

\newcommand{\bcl}{\begin{cl}}
\newcommand{\ecl}{\end{cl}}

\newcommand{\bmlem}{\begin{mlem}}
\newcommand{\emlem}{\end{mlem}}

\newcommand{\bscl}{\begin{scl}}
\newcommand{\escl}{\end{scl}}

\newcommand{\bcons}{\begin{conjs}}
\newcommand{\econs}{\end{conjs}}

\newcommand{\bprop}{\begin{prop}}
\newcommand{\eprop}{\end{prop}}

\newcommand{\br}{\begin{rem}}
\newcommand{\er}{\end{rem}}
\newcommand{\brs}{\begin{rems}}
\newcommand{\ers}{\end{rems}}
\newcommand{\bo}{\begin{obser}}
\newcommand{\eo}{\end{obser}}
\newcommand{\bos}{\begin{obsers}}
\newcommand{\eos}{\end{obsers}}
\newcommand{\bpf}{\begin{pf}}
\newcommand{\epf}{\end{pf}}
\newcommand{\ba}{\begin{array}}
\newcommand{\ea}{\end{array}}
\newcommand{\beq}{\begin{eqnarray}}
\newcommand{\beqq}{\begin{eqnarray*}}
\newcommand{\eeq}{\end{eqnarray}}
\newcommand{\eeqq}{\end{eqnarray*}}

\newcounter{minutes}
\divide\time by 60
\newcounter{hours}
\multiply\time by 60 \addtocounter{minutes}{-\time}

\begin{document}

\bibliographystyle{amsplain}

\title{$L^p\to L^q$ norm estimates of Cauchy transforms on the Dirichlet problem and their applications}
\author[J.-F. Zhu]{Jian-Feng Zhu}
\address{Jian-Feng Zhu, School of Mathematical Sciences, Huaqiao
University, Quanzhou 362021, China and Department of Mathematics,
Shantou University, Shantou, Guangdong 515063, China.}
\email{flandy@hqu.edu.cn}

\author[A. Rasila]{Antti Rasila}
\address{Antti Rasila, Technion -- Israel Institute of Technology, Guangdong Technion, Shantou, Guangdong 515063, China.}
\email{antti.rasila@iki.fi; antti.rasila@gtiit.edu.cn}

\date{\today}

\begin{abstract}
Denote by $C^{\alpha}(\ID)$ the space of the functions $f$ on the unit disk $\ID$ which are H\"older continuous with the exponent $\alpha$, and
denote by $C^{1, \alpha}(\ID)$ the space which consists of differentiable functions $f$ such that their derivatives are in the space $C^{\alpha}(\ID)$.
Let $\mathcal{C}$ be the Cauchy transform of Dirichlet problem.
In this paper, we obtain the norm estimates of $\|\mathcal{C}\|_{L^p\to L^q}$, where
$3/2<p<2$ and $q=p/(p-1)$. As an application, we show that
if $3/2<p<2$, then $u\in C^{\mu}(\ID)$, where $\mu=2/p-1$. We also show that
if $2<p<\infty$, then $u\in C^{1, \nu}(\ID)$, where $\nu=1-2/p$.
Finally, for the case $p=\infty$, we show that $u$ is not necessarily in $C^{1, 1}(\ID)$, but its gradient, i.e., $|\nabla u|$
is Lipschitz continuous with respect to the pseudo-hyperbolic metric.
This paper is inspired by \cite[Chapter 4]{Astala} and \cite{kalaj1}.
\end{abstract}

\maketitle 
\pagestyle{myheadings} \markboth{J.-F. Zhu and A. Rasila}{$L^p\to L^q$ norm estimates of Cauchy transforms on the Dirichlet problem and their applications}


\section{Introduction}\label{sec-1}
\subsection*{The space $L^{p} (\Omega)$}
Throughout this paper, we use $\ID$ denote by the unit disk, and $\mathbb{T}$ the unit circle.
Suppose $\Omega$ is a domain in the complex plane $\mathbb{C}$.
Denote by $L^p(\Omega)$ the space of the
complex-valued measurable functions on $\Omega$ with a finite integral
$$\|f\|_{p}=\left(\int_{\Omega}|f(z)|^p\mathrm{d} A(z)\right)^{\frac{1}{p}},\ \ \ 1\leq p<\infty,$$
where $\mathrm{d}A(z)$
is the normalized area measure on $\Omega$ (cf. \cite[Page 1]{Hendenmalm}). For the case $p=\infty$, we let $L^{\infty}(\Omega)$ denote the space of (essentially) bounded
functions on $\Omega$. For $f\in L^{\infty}(\Omega)$, we define
$$\|f\|_\infty=\esssup\{|f(z)|: z\in\Omega\}.$$
It is known that the space $L^{\infty}(\ID)$ is a Banach space with the above norm  (cf. \cite[Page 2]{Hendenmalm}).

\subsection*{The Sobolev space $W^{1, p}(\Omega)$}
For $k\geq0$ and $p\geq1$, the Sobolev space $W^{1, p}(\Omega)$ is the Banach space of $k$-times weak differentiable $p$-integrable functions. The norm in $W^{1, p}(\Omega)$
is defined by
$$\|u\|_{W^{1, p}}=\left(\int_{\Omega}\sum_{|\alpha|\leq k}|D^{\alpha}u|^p\mathrm{d}A\right)^{1/p}.$$
The space $W_0^{1, p}(\Omega)$ is obtained by taking the closure of $C^k_0(\Omega)$ in $W^{k, p}(\Omega)$, where
$C_0^{k, p}(\Omega)$ is the space of $k$ times continuously differentiable functions with a compact support in $\Omega$ (cf. \cite[Pages 153--154]{GT}).

\subsection*{The Cauchy transform of a solution to the Dirichlet problem}
The Poisson equation is given as follows:
\be\label{poisson-equation}
\left\{
\begin{array}
{r@{\ }l}
&u_{z\bar{z}}=g(z), \ \   z\in\Omega,\\
\\
&u\in W^{1, p}_0,
\end{array}\right.\ee
where $\Delta u=4u_{z\bar{z}}$ is the Laplacian of $u$ and $g\in L^p(\Omega)$.
It is known that if $\Omega=\ID$ and $g\in L^p(\ID)$, where $1\leq p\leq\infty$, then  the weak solution of the Poisson equation is
$$u(z)=G[g](z)=\int_{\ID}\log\left|\frac{z-\tau}{1-\bar{z}\tau}\right|^2g(\tau)\mathrm{d}A(\tau),$$
where $z$, $\tau\in \ID$, and
$$G(z, \tau)=\log\left|\frac{z-\tau}{1-\bar{z}\tau}\right|^2$$
is the Green function.

Suppose $u=G[g]$ is a solution to (\ref{poisson-equation}), where $g\in L^p(\ID)$ and $1<p<\infty$. We may define the
{\it Cauchy transform of the solution to the Dirichlet problem} as follows:
$$\mathcal{C}[g](z)=\frac{\partial}{\partial z}u=\int_{\ID}\left(\frac{1}{z-\tau}+\frac{\bar{\tau}}{1-z\bar{\tau}}\right)g(\tau)\mathrm{d}A(\tau).$$
The operator $\mathcal{C}$ is then induced by the complex partial $z$-derivative of the Green's function (cf. \cite[Page 155]{Astala}).

Let $\mathfrak{C}$ be the {\it Cauchy integral operator} ({\it the Cauchy transform}) which is defined as follows:
$$\mathfrak{C}[g](z)=\int_{\Omega}\frac{g(\tau)}{\tau-z}\mathrm{d}A(\tau).$$
The following integral operator $\mathfrak{J}_0^*$ was introduced in \cite[Page 12]{Baranov}  and is given as follows:
$$\mathfrak{J}_0^*[g](z)=\int_{\ID}\frac{\bar{\tau}}{1-z\bar{\tau}}g(\tau)\mathrm{d}A(\tau).$$
Now, it is easy to see that $\mathcal{C}=\mathfrak{J}_0^*-\mathfrak{C}$ (cf. \cite{Astala, Baranov, kalaj1}).
Moreover, elementary calculations show that (cf. \cite[(1.4) and (1.5)]{kalaj1})
$$\mathcal{C}[g](z)=\frac{\partial}{\partial z}u=\int_{\ID}\frac{1-|\tau|^2}{(z-\tau)(1-z\bar{\tau})}g(\tau)\mathrm{d}A(\tau),$$
and
$$\overline{\mathcal{C}}[g](z)=\frac{\partial}{\partial \bar{z}}u=\int_{\ID}\frac{1-|\tau|^2}{(\bar{z}-\bar{\tau})(1-\bar{z}\tau)}g(\tau)\mathrm{d}A(\tau).$$

Recall that the standard operator norm of an operator $T: X\rightarrow Y$ between normed spaces $X$ and $Y$ is defined by
$$\|T\|_{X\rightarrow Y}=\sup\{\|Tx\|_Y:\|x\|_X=1\}.$$
For simplicity, if $X=Y=L^{p}(\ID)$, then we write $\|T\|_p$ instead of $\|T\|_{L^{p}\rightarrow L^{p}}$ for
the $L^p$ norm of the operator $T$.

For $p>2$, it was shown in \cite[Page 155]{Astala} that $\mathcal{C}[g](z)$ is a continuous function on the closed disk $\overline{\ID}$.
In 2012, Kalaj proved in \cite[Theorem B]{kalaj1} that there exists a constant $M_p$ depending only on $p$,
such that $\|\mathcal{C}\|_{L^p(\ID)\to L^\infty(\ID)}=M_p$.
Moreover, he obtained norm estimates for $\|\mathcal{C}\|_p$, and showed that the results are sharp for $p=1, 2$ and $\infty$ (cf. \cite[Theorem A]{kalaj1}).

\subsection*{The pseudo-hyperbolic distance on $\ID$}
For each $z\in\ID$, let $\varphi_w$ denote the M\"{o}bius transform of the form
$$\varphi_w(z)=\frac{w-z}{1-\bar{w}z},$$
where $w\in\ID$. The {\it pseudo-hyperbolic distance} on $\ID$ is defined as follows
(cf. \cite{Proceeding-AMS}):
\be\label{hz-2018-eq-1}
\rho(z, w)=|\varphi_w(z)|, \ \ \ z, w\in\ID.
\ee
We note that the pseudo-hyperbolic distance is invariant under M\"{o}bius transformations, that is,
$$\rho(f(z), f(w))=\rho(z, w),$$
for any $f\in$Aut$(\ID)$, the M\"obius automorphism of $\ID$, where $z, w\in\ID$ (cf. \cite{Proceeding-AMS}). Moreover, it has the following useful property:
\be\label{hz-2018-eq-2}
1-\rho(z, w)^2=\frac{(1-|z|^2)(1-|w|^2)}{|1-\bar{z}w|^2}.
\ee

\begin{defn}$($cf. \cite[Page 115]{Astala} $)$
The H\"older spaces $C^{\mu}(\mathbb{C})$, $0<\mu\leq1$, consist of continuous functions $f:\mathbb{C}\to\mathbb{C}$ that satisfy the H\"older
condition
$$\|f\|_{C^{\mu}(\mathbb{C})}=\sup\limits_{z\neq w}\frac{|f(z)-f(w)|}{|z-w|^{\mu}}<\infty.$$
The space $C^{1, \nu}(\mathbb{C})$, $0<\nu\leq1$, consists of $C^1$ functions $f:\mathbb{C}\to\mathbb{C}$ that satisfy the following condition:
$$\|f\|_{C^{1, \nu}(\mathbb{C})}=\sup\limits_{z\neq w}\frac{|\nabla f(z)-\nabla f(w)|}{|z-w|^{\nu}}<\infty.$$
\end{defn}

\subsection*{Main results}
In this paper, we show that if $g\in L^p(\ID)$, where $3/2<p<2$, then the operator $\mathcal{C}$ will map $L^p(\ID)$ to $L^q(\ID)$, where $q=p/(p-1)$
is the conjugate exponent of $p$. This result partly improves the corresponding results in \cite[Theorem A]{kalaj1}. As an application, by using Morrey's inequality (cf. \cite{Ryan, Morrey}), we have $G[g]\in C^{\mu}(\ID)$ is H\"older continuous with the exponent $\mu=1-2/q$. Moreover, we prove that if $2<p<\infty$ and $g\in L^p(\ID)$, then $G[g]\in C^{1, \nu}(\ID)$, where $\nu=1-2/p$.
For the case $p=\infty$, we give Example \ref{example-2} which shows that $\mathcal{C}$ is not of $C^1(\ID)$ space or even
the space of Lipschitz continuous functions. This implies that when $g\in L^{\infty}(\ID)$, its Green potential $G[g]$ is not in the space $C^{1, 1}(\ID)$,
or even in $C^{1, \textbf{Lip}}(\ID)$, the space of the functions with Lipschitz continuous derivatives. However, by applying the pseudo-hyperbolic distance, we show that
$|\nabla G[g]|$ is Lipschitz continuious with respect to the pseudo-hyperbolic metric.

More precisely, our results are as follows:

\begin{thm}\label{20-March-10-thm-1}
Suppose $g\in L^p(\ID)$, where $3/2< p<2$, and $q=p/(p-1)$. Then $\mathcal{C}$ is an operator of $L^p(\ID)$ to $L^q(\ID)$.
Moreover, we have
$$\|\mathcal{C}\|_{L^p\to L^q}^q\leq 2^{3/2}\left(\frac{2}{2-q/p}+\frac{1}{\Gamma(2-q/p)}\right).$$
\end{thm}

The following example shows that for the case $p_0=3/2$ and $q_0=3$, there exists a function $u=G[g]$, such that $g\in L^{p_0}(\ID)$, but $\mathcal{C}[g](z)=\frac{\partial}{\partial z}u(z)\notin L^{q_0}(\ID)$.

\begin{figure}
\centering
\includegraphics[width=4.1cm]{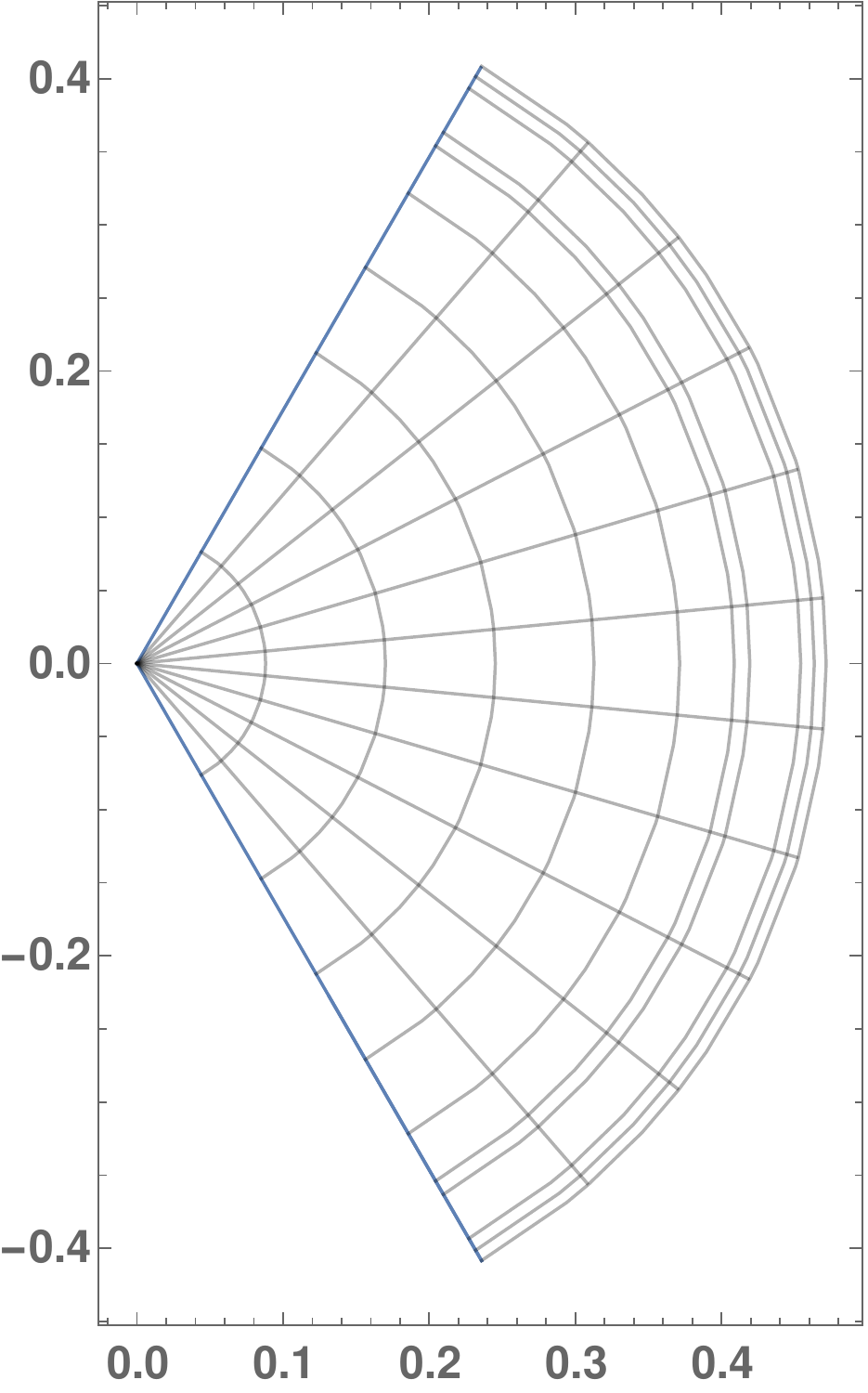}
\quad
\includegraphics[width=4cm]{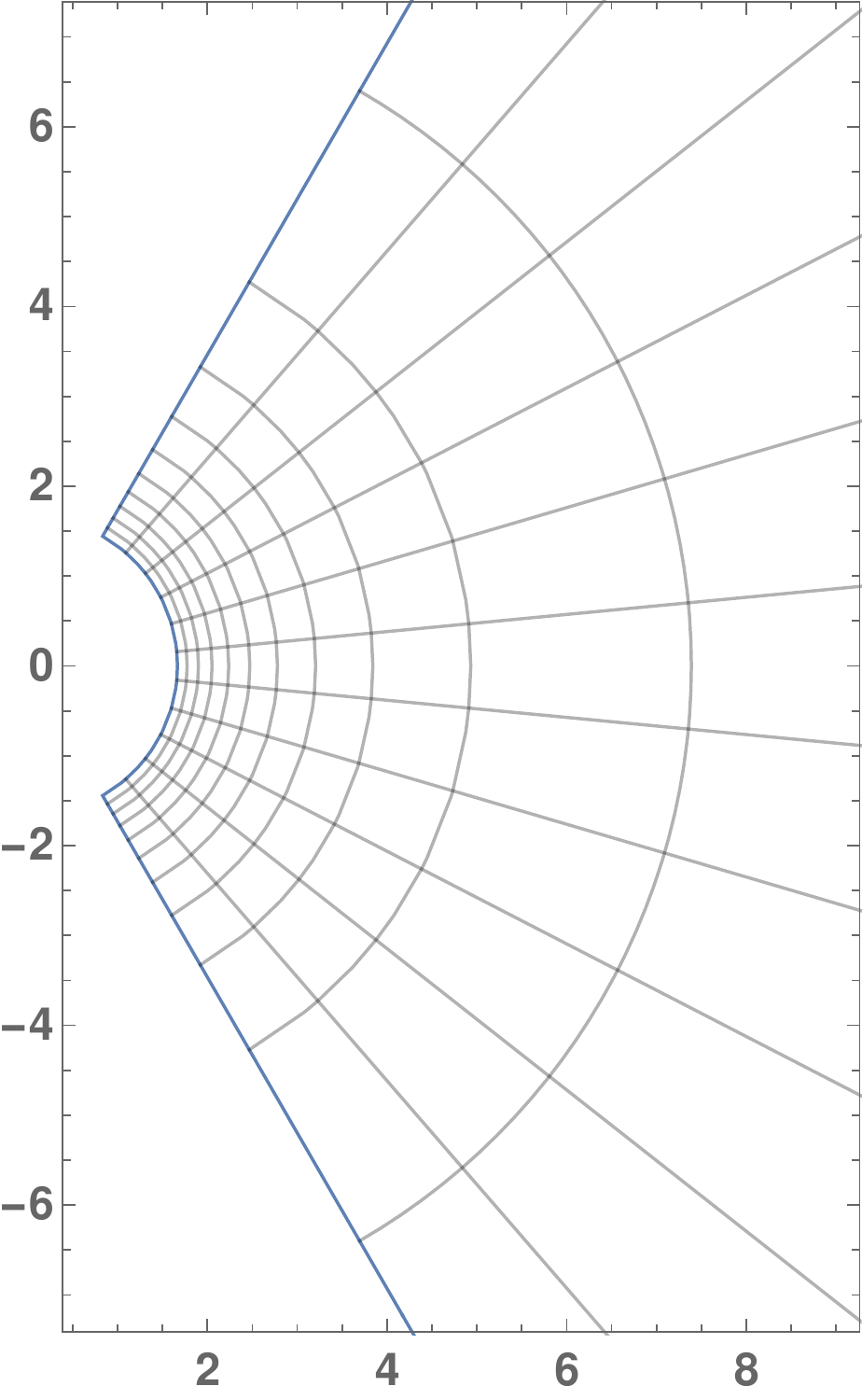}
\quad
\includegraphics[width=4.4cm]{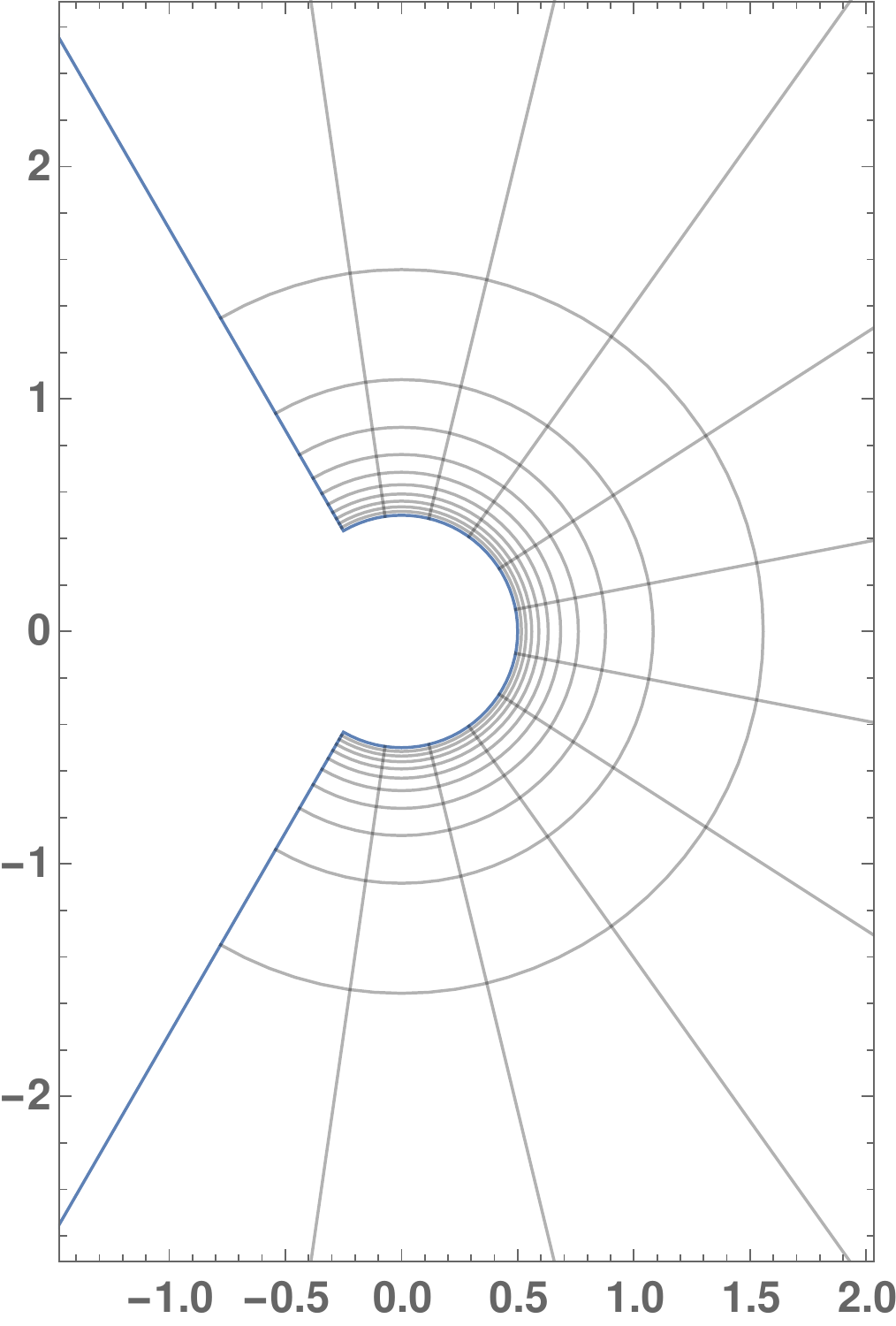}

\caption{The images of $\mathbb{D}$ under the mappings $u(z),g(z)$, and $u_z(z)$ of Example \ref{ex1}. }\label{ex1fig}
\end{figure}

\begin{exam}\label{ex1}
Let $u(z)=z^{\frac{1}{3}}(1-|z|^{2\alpha})$, where $z\in\ID$ and $\alpha>\frac{1}{6}$. Then $u(\zeta)=0$, for any $\zeta=e^{i\theta}\in\mathbb{T}$, and
$g(z)=\Delta u=-4\alpha(\alpha+\frac{1}{3})z^{\frac{1}{3}}|z|^{2\alpha-2}$. Elementary calculations show that $g\in L^{\frac{3}{2}}(\ID)$ and
$\frac{\partial}{\partial z}u(z)=z^{-\frac{2}{3}}(\frac{1}{3}-(\frac{1}{3}+\alpha)|z|^{2\alpha})\notin L^3(\ID)$. See Figure \ref{ex1fig}.
\end{exam}

As an application of Theorem \ref{20-March-10-thm-1}, we have the following results:
\begin{thm}\label{20-March-thm-1}
Suppose $g\in L^p(\ID)$, where $3/2<p<\infty$, and $G[g]$ is the Green potential of $g$.
Then\emph{:}
\begin{itemize}
  \item [$(1)$] for $3/2<p<2$, $G[g]\in C^{\mu}(\ID)$, where $\mu=2/p-1$;
  \item [$(2)$] for $2<p<\infty$, $G[g]\in C^{1, \nu}(\ID)$, where $\nu=1-2/p$.
\end{itemize}
\end{thm}

For the case $p=\infty$, the following example shows that $G[g]$ does not necessarily belong to the space $C^{1, 1}(\ID)$ or even to $C^{1, \textbf{Lip}}(\ID)$.

\begin{figure}
\centering
\includegraphics[width=7cm]{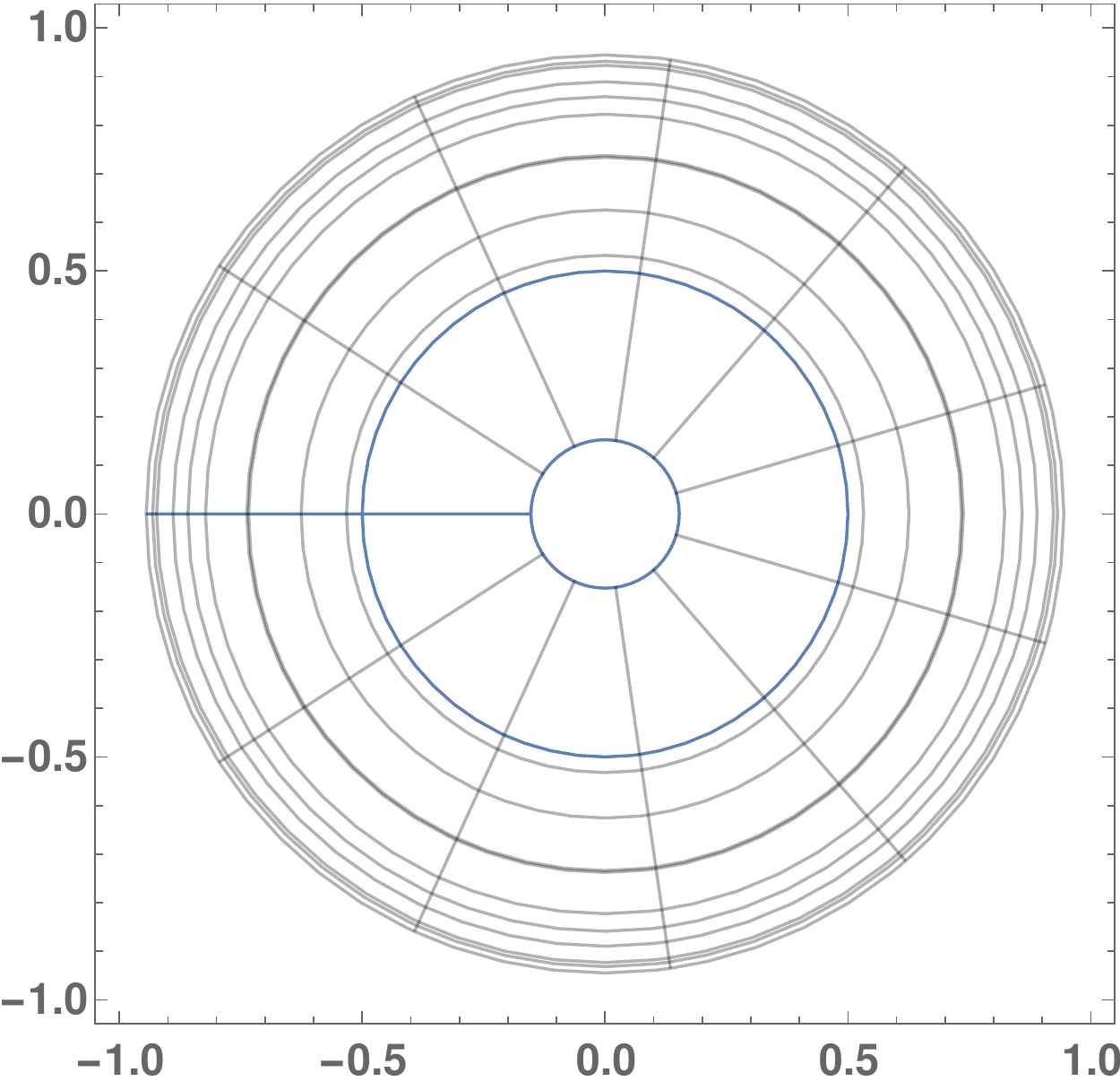}

\caption{The image of $\mathbb{D}$ under the mapping $\mathcal{C}(z)$ of Example \ref{example-2}. }\label{ex2fig}
\end{figure}

\begin{exam}\label{example-2} (cf. \cite[Page 116]{Astala})
Let $g(z)=z/\bar{z} \in L^{\infty}(\ID)$. Then $\mathfrak{C}[g](z)=-z\log|z|^2$ and $\mathfrak{J}_0^*[g](z)=z/2$.
This shows that $\mathcal{C}=\mathfrak{J}_0^*-\mathfrak{C}$ does not take $L^{\infty}(\ID)$ to the space $C^1(\ID)$ or even
the space of Lipschitz continuous functions. See Figure \ref{ex2fig}.
\end{exam}

Next, we show that $|\nabla G[g]|$ is Lipschitz continuous with respect to the pseudo-hyperbolic metric.
Suppose that $g\in L^{\infty}(\ID)$ and that $u(z)=G[g](z)$ is the Green potential of $g$. Then the gradient of $u$ is defined by
$$\Lambda_u(z)=|u_z(z)|+|u_{\bar{z}}(z)|=\big|\mathcal{C}[g](z)\big|+\big|\overline{\mathcal{C}}[g](z)\big|.$$
We have the following result:

\begin{thm}\label{thm-3-18}
Let $u=G[g]$ be the Green potential of $g$, where $g\in L^{\infty}(\ID)$. Then
\be\label{DZ-thm-3-eq-1}
\big|(1-|z|^2)\Lambda_u(z)-(1-|w|^2)\Lambda_u(w)\big|\leq2\|g\|_{\infty}\rho(z, w),
\ee
holds for all $z, w\in\ID$.
\end{thm}

The rest of this paper is organized as follows: In Section \ref{sec-2}, we recall some known results and prove two lemmas which will be used in proving Theorem \ref{20-March-10-thm-1}. In Section \ref{sec-2}, we prove Theorems \ref{20-March-10-thm-1} to \ref{thm-3-18}. In Section \ref{sec-3}, we provide some additional insights related
to the constant $C_q$, which occurs in Theorem \ref{20-March-thm-1}.
\section{Auxiliary results}\label{sec-2}
In this section, we recall some known results and prove two lemmas which will be used in proving our main theorems.
We start with the following Morrey's inequality (cf. \cite{Ryan, Morrey}).

\begin{Thm}\label{Morrey-ineq}$($Morrey's inequality$)$
Assume that $n<p<\infty$, and let $u\in D^{1, p}(\mathbb{R}^n)$ $($i.e., the derivative of $u$ exists and of $L^p(\mathbb{R}^n)$ space$)$.
Then there exists a constant $C>0$ such that
$$\sup\limits_{x\neq y}\left\{\frac{|u(x)-u(y)|}{|x-y|^{1-n/p}}\right\}\leq C\left(\int_{\mathbb{R}^n}|Du|^p\mathrm{d}x\right)^{1/p}.$$
\end{Thm}

The following equality (\ref{bz-lem-2.0}) and Lemma \Ref{DZ-2019-lem-0} will be applied in the proofs of Lemmas \ref{lem-1-20-March-10}, \ref{lem-2-20-March-10}, and
Theorem \ref{20-March-10-thm-1}.

For $\beta>0$, $z\in\mathbb{D}$ and $\zeta=e^{i\theta}\in\mathbb{T}$, we have
$$\frac{1}{(1-z\zeta)^\beta}=\sum_{n=0}^\infty\frac{\Gamma(n+\beta)}{n!\Gamma(\beta)}z^n\zeta^n.$$
By using {\it Parseval's theorem}, one gets the identity:
\be\label{bz-lem-2.0}\frac{1}{2\pi}\int_0^{2\pi}\frac{\mathrm{d}\theta}{|1-ze^{i\theta}|^{2\beta}}=\sum_{n=0}^\infty\left(\frac{\Gamma(n+\beta)}{n!\Gamma(\beta)}\right)^2|z|^{2n}.\ee

Recall the following estimate:

\begin{Lem}$($\emph{cf.} \cite[(7)]{Gautschi}$)$\label{DZ-2019-lem-0}
For $0\leq a\leq1$ and $n=1, 2, \ldots$,
$$\frac{1}{(n+1)^{1-a}}\leq\frac{\Gamma(n+a)}{n!}\leq\frac{1}{n^{1-a}}.$$
\end{Lem}

Next, we prove two lemmas which will be used in proving our Theorem \ref{20-March-10-thm-1}.
\begin{lem}\label{lem-1-20-March-10}
For $1\leq\beta<2$, let
\be\label{Ibeta-def}I_{\beta}=\int_{\ID}\left(\frac{1-|w|^2}{|z-w||1-\bar{z}w|}\right)^{\beta}\mathrm{d}A(w),\ee
where $z\in\ID$. Then
\be\label{Ibeta}I_{\beta}\leq2^{1+\beta/2}\Gamma(1+\beta)\Gamma(2-\beta).\ee
\end{lem}
\bpf
By using the M\"obius transformation $\eta=\frac{z-w}{1-\bar{z}w}$, we obtain
$$I_\beta=(1-|z|^2)^{2-\beta}\int_{\ID}\left(\frac{1-|\eta|^2}{|\eta|}\right)^{\beta}\frac{1}{|1-\bar{z}\eta|^4}\mathrm{d}A(\eta).$$
Suppose that $\eta=re^{it}\in\ID$. By applying (\ref{bz-lem-2.0}) and the following equality:
$$\int_0^1r^{2n+1-\beta}(1-r^2)^\beta \mathrm{d}r=\frac{\Gamma(1+\beta)\Gamma(n+1-\beta/2)}{2\Gamma(n+2+\beta/2)},\ \ \ n=0, 1, \ldots,$$
we have
\begin{align}\label{March-12}
  I_{\beta} & = \frac{(1-|z|^2)^{2-\beta}}{\pi}\int_{0}^1r^{1-\beta}(1-r^2)^{\beta}\mathrm{d}r\int_0^{2\pi}\frac{1}{|1-\bar{z}re^{it}|^4}\mathrm{d}\theta\\ \nonumber
   & =(1-|z|^2)^{2-\beta}\Gamma(1+\beta)\sum\limits_{n=0}^{\infty}(n+1)^2\frac{\Gamma(n+1-\beta/2)}{\Gamma(n+2+\beta/2)}|z|^{2n}.
\end{align}
Note that for $0<a<1$, the formula
\be\label{eq-series-10}\frac{1}{(1-z)^a}=\sum\limits_{n=0}^{\infty}\frac{\Gamma\left(n+a\right)}{n!\Gamma\left(a\right)}z^n\ee
holds for every $z\in\ID$.
Moreover, according to Lemma \Ref{DZ-2019-lem-0}, one obtains the following inequality:
\be\label{eq-March-10-1}(n+1)^2\frac{\Gamma(n+1-\beta/2)}{\Gamma(n+2+\beta/2)}\leq2^{1+\beta/2}\frac{\Gamma(n+2-\beta)}{n!},\ \ \ n=1, 2, \ldots.\ee
It is easy to see that when $n=0$, the above inequality (\ref{eq-March-10-1}) still holds.
Combining (\ref{March-12}), (\ref{eq-series-10}) and (\ref{eq-March-10-1}), we see that (\ref{Ibeta}) holds, which completing the proof.
\epf

\begin{lem}\label{lem-2-20-March-10}
For $1\leq\beta<2$, let
$$J_{\beta}=\int_{\ID}\left(\frac{1-|w|^2}{|z-w||1-\bar{z}w|}\right)^{\beta}\mathrm{d}A(z),$$
where $w\in\ID$. Then
\be\label{Jbeta}J_{\beta}\leq\frac{2}{2-\beta}+\frac{1}{\Gamma(2-\beta)}.\ee
\end{lem}
\bpf
Following the proof of Lemma \ref{lem-1-20-March-10}, by letting $z=\frac{w-\eta}{1-\bar{w}\eta}$ and applying (\ref{bz-lem-2.0}), one has
\begin{align*}
  J_\beta & =(1-|w|^2)^{2-\beta}\int_{\ID}\frac{1}{|\eta|^{\beta}}\frac{1}{|1-\bar{z}\eta|^{4-2\beta}}\mathrm{d}A(\eta) \\
   &=2(1-|w|^2)^{2-\beta}\sum\limits_{n=0}^{\infty}\left(\frac{\Gamma(n+2-\beta)}{n!\Gamma(2-\beta)}\right)^2\frac{|w|^{2n}}{2n+2-\beta}\\
   &=2(1-|w|^2)^{2-\beta}\left(\frac{1}{2-\beta}+\sum\limits_{n=1}^{\infty}\left(\frac{\Gamma(n+2-\beta)}{n!\Gamma(2-\beta)}\right)^2\frac{|w|^{2n}}{2n+2-\beta}\right).
\end{align*}
Since $1\leq\beta<2$, by using Lemma \Ref{DZ-2019-lem-0} we see that
$$\frac{\Gamma(n+2-\beta)}{n!(2n+2-\beta)}<\frac{1}{2n^\beta}\leq\frac{1}{2},\ \ \ n=1, 2, \ldots.$$
Thus
$$J_\beta\leq2(1-|w|^2)^{2-\beta}\left(\frac{1}{2-\beta}+\frac{1}{2\Gamma(2-\beta)}\sum\limits_{n=1}^{\infty}\frac{\Gamma(n+2-\beta)}{n!\Gamma(2-\beta)}|w|^{2n}\right).$$
Again, by (\ref{eq-series-10}), we have
\begin{align*}
  J_\beta & \leq2(1-|w|^2)^{2-\beta}\left[\frac{1}{2-\beta}+\frac{1}{2\Gamma(2-\beta)}\left(\frac{1}{(1-|w|^2)^{2-\beta}}-1\right)\right] \\
   & \leq \frac{2}{2-\beta}+\frac{1}{\Gamma(2-\beta)}.
\end{align*}
Therefore, the desired inequality (\ref{Jbeta}) follows.
\epf

\section{Proofs of the main results}\label{sec-3}
\subsection*{Proof of Theorem \ref{20-March-10-thm-1}}
Recall that
$$\mathcal{C}[g](z)=\int_{\ID}\frac{1-|\tau|^2}{(z-\tau)(1-z\bar{\tau})}g(\tau)\mathrm{d}A(\tau).$$
By using Lemma \ref{lem-1-20-March-10} and the H\"older's inequality for integrals, we have
$$|\mathcal{C}[g](z)|\leq I_1^{1/q}\left(\int_{\ID}\frac{1-|\tau|^2}{|z-\tau||1-\bar{z}\tau|}|g(\tau)|^p\mathrm{d}A(\tau)\right)^{1/p},$$
where $q=\frac{p}{p-1}$. The assumption $g\in L^p(\ID)$ ensures that
$$\int_{\ID}|g(\tau)|^p\mathrm{d}A(\tau)=\|g\|_p^p<\infty,$$
and thus,
$$|g(\tau)|^p\frac{\mathrm{d}A(\tau)}{\|g\|_p^p}$$
is a probability measure in $\ID$.

Observe that under the assumption of $3/2<p<2$, we have $1<q/p<2$.
Applying Jensen's inequality (cf. \cite[Page 231]{kalaj1}), we obtain
\begin{align}\label{March-12-2}
  |\mathcal{C}[g](z)|^q & \leq I_1\,\|g\|_{p}^{q}\left(\int_{\ID}\frac{1-|\tau|^2}{|z-\tau||1-\bar{z}\tau|}|g(\tau)|^p\frac{\mathrm{d}A(\tau)}{\|g\|_p^p}\right)^{q/p}\\ \nonumber
   & \leq I_1\,\|g\|_{p}^{q-1}\int_{\ID}\left(\frac{1-|\tau|^2}{|z-\tau||1-\bar{z}\tau|}\right)^{q/p}|g(\tau)|^p\mathrm{d}A(\tau).
\end{align}
It follows from Lemma \ref{lem-1-20-March-10} and the assumption $g\in L^p(\ID)$ that
$$\left(\frac{1-|\tau|^2}{|z-\tau||1-\bar{z}\tau|}\right)^{q/p}|g(\tau)|^p\in L^1(\ID\times \ID).$$
By using Lemma \ref{lem-2-20-March-10} and (\ref{March-12-2}), one obtains
\begin{align*}
  \int_{\ID}|\mathcal{C}[g](z)|^q\mathrm{d}A(z) & \leq I_1\,\|g\|_{p}^{q-1}\left(\int_{\ID}|g(\tau)|^p\mathrm{d}A(\tau)\right)\int_{\ID}\left(\frac{1-|\tau|^2}{|z-\tau||1-\bar{z}\tau|}\right)^{q/p}\mathrm{d}A(z)\\
   & \leq \|g\|_p^{q-1+p}I_1J_{q/p}<\infty.
\end{align*}
This shows that $\mathcal{C}[g]\in L^q(\ID)$, where $q=p/(p-1)$.
Moreover, we have
\begin{align*}
  \|\mathcal{C}\|_{L^p\to L^q}^q & =\sup\{\|\mathcal{C}[g]\|_q^q: \|g\|_p=1\} \\
   & \leq  2^{3/2}\left(\frac{2}{2-q/p}+\frac{1}{\Gamma(2-q/p)}\right).
\end{align*}
The proof is complete. \qed
\subsection*{Proof of Theorem \ref{20-March-thm-1}}
(1) The assumption $3/2<p<2$ ensures that $2<q=p/(p-1)<2p$. According to Theorem  \ref{20-March-10-thm-1}, we see that
$\mathcal{C}[g]\in L^q(\ID)$. Similarly, we may prove that $\overline{\mathcal{C}}[g]\in L^q(\ID)$.
Then, by Theorem \Ref{Morrey-ineq}, we see that $G[g]\in C^{\mu}(\ID)$ is H\"older continuous with the exponent $\mu=1-2/q=2/p-1$.

(2) Suppose $g\in L^p(\ID)$, where $2<p<\infty$. Then, for $z$, $w\in\ID$ satisfying $z\neq w$, we have
\begin{align*}
  \mathcal{C}[g](z)-\mathcal{C}[g](w) & =\mathfrak{J}_0^*[g](z)-\mathfrak{J}_0^*[g](w)-\big[\mathfrak{C}[g](z)-\mathfrak{C}[g](w)\big] \\
   & =\int_{\ID}\left[\frac{(z-w)\bar{\tau}^2}{(1-z\bar{\tau})(1-w\bar{\tau})}+\frac{w-z}{(z-\tau)(w-\tau)}\right]g(\tau)\mathrm{d}A(\tau).
\end{align*}

We first estimate $\mathfrak{J}_0^*[g](z)-\mathfrak{J}_0^*[g](w)$ as follows:
\begin{align*}
  |\mathfrak{J}_0^*[g](z)-\mathfrak{J}_0^*[g](w)| & \leq|z-w|\int_{\ID}\frac{|g(\tau)|}{|\bar{z}-1/\tau||\bar{w}-1/\tau|}\mathrm{d}A(\tau) \\
   &\leq |z-w|\int_{\mathbb{C}\setminus \ID}\frac{|g(1/\eta)|}{|\bar{z}-\eta||\bar{w}-\eta|}\mathrm{d}A(\eta),
\end{align*}
where the last inequality holds because $1/{|\eta|^4}\leq1$, for any $\eta\in \mathbb{C}\setminus \ID$.
Because $g\in L^p(\ID)$ ($2<p<\infty$), by using H\"older inequality for integrals, we have
\begin{align*}
  |\mathfrak{J}_0^*[g](z)-\mathfrak{J}_0^*[g](w)| & \leq\|g\|_p|z-w|\left(\int_{\mathbb{C}\setminus\ID}\frac{\mathrm{d}A(\eta)}{|\bar{z}-\eta|^q|\bar{w}-\eta|^q}\right)^{1/q} \\
  & \leq|z-w|^{2/q-1}\|g\|_p\left(\int_{\mathbb{C}}\frac{\mathrm{d}A(\xi)}{|\xi|^q|1-\xi|^q}\right)^{1/q}\\
  &=C_q\|g\|_p|z-w|^{1-2/p},
\end{align*}
where $q=p/(p-1)\in(1, 2)$, and
$$C_q=\left(\int_{\mathbb{C}}\frac{\mathrm{d}A(\xi)}{|\xi|^q|1-\xi|^q}\right)^{1/q}$$
is a constant depending only on $q$ (see Section \ref{sec-app} for more details).

Next, we estimate $\mathfrak{C}[g](z)-\mathfrak{C}[g](w)$ as follows (cf. \cite[Theorem 4.3.13]{Astala}):
\begin{align*}
  |\mathfrak{C}[g](z)-\mathfrak{C}[g](w)| & \leq\|g\|_p|z-w|\left(\int_{\ID}\frac{1}{|z-\tau|^q|w-\tau|^q}\mathrm{d}A(\tau)\right)^{1/q} \\
   & \leq\|g\|_p|z-w|^{2/q-1}\left(\int_{\mathbb{C}}\frac{1}{|\eta|^q|1-\eta|^q}\mathrm{d}A(\eta)\right)^{1/q}\\
   &=C_q\|g\|_p|z-w|^{1-2/p}.
\end{align*}
Then
$$|\mathcal{C}[g](z)-\mathcal{C}[g](w)|\leq2C_q\|g\|_p|z-w|^{1-2/p}.$$
Similarly, we obtain the following estimate:
$$|\overline{\mathcal{C}}[g](z)-\overline{\mathcal{C}}[g](w)|\leq2C_q\|g\|_p|z-w|^{1-2/p}.$$
Because $\mathcal{C}[g](z)=\frac{\partial}{\partial z}G[g](z)$ and $\overline{\mathcal{C}}[g](z)=\frac{\partial}{\partial \bar{z}}G[g](z)$, we see that $G[g]\in C^{1, \nu}(\ID)$, where $\nu=1-2/p$.
The proof of Theorem \ref{20-March-thm-1} is complete.
\qed

\subsection*{Proof of Theorem \ref{thm-3-18}}
For every $z,w\in \ID$, set $\lambda =\varphi_w(z)=(w-z)/(1-\overline{w}z)$ and let $\psi=u\circ \varphi _w$. Then
\be\label{phi-3-18-1}(1-|w|^2)\Lambda _u(w)=|\psi_z(0)|+|\psi_{\bar{z}}(0)|=\Lambda _{\psi }(0).\ee
On the other hand, as $z=(w-\lambda)/(1-\bar{w}\lambda)=\varphi_w(\lambda)$, it follows from (\ref{hz-2018-eq-2}) that
\begin{eqnarray*}
(1-|z|^2)\Lambda _u(z)&=&(1-|\varphi_w(\lambda )|^2)\Lambda _u(\varphi _w(\lambda ))\\ \nonumber
&=&(1-|\lambda |^2 )|\varphi' _w(\lambda )|\Lambda _u(\varphi _w(\lambda ))\\ \nonumber
&=&(1-|\lambda |^2)\Lambda _{\psi }(\lambda ).
\end{eqnarray*}
Therefore,
\begin{align}\label{main-1}
\nonumber \left|(1-|z|^2)\Lambda _u(z)-(1-|w|^2)\Lambda _u(w)\right|&=\left|(1-|\lambda |^2)\Lambda _{\psi }(\lambda )-\Lambda _{\psi }(0)\right|\\
&\leq |\lambda |^2\Lambda _{\psi }(0)+(1-|\lambda |^2)\left|\Lambda _{\psi}(\lambda )-\Lambda _{\psi}(0)\right|.
\end{align}
It follows from \cite[Lemma 2.7]{DKal2011TAMS} that for any $g\in L^{\infty}(\ID)$, we have
$\Lambda_{u}(z)\leq 2/3\|g\|_{\infty}.$ By using (\ref{phi-3-18-1}), we see that
\be\label{eq-march-18-1}\Lambda_{\psi}(0)\leq\frac{2}{3}\|g\|_{\infty}.\ee

Next, we estimate $\left|\Lambda _{\psi}(\lambda )-\Lambda _{\psi}(0)\right|$ as follows:
Let $\Gamma$ be the line segment joining $0$ and $\lambda$, i.e., $\Gamma$ has the parametric equation $z(t)=t\lambda$, where $0\leq t\leq1$. Then
\begin{align}\label{march-18-2}
  \left|\Lambda _{\psi}(\lambda )-\Lambda _{\psi}(0)\right| & \leq\int_{\Gamma}\Lambda_u(\varphi_w(t\lambda))|\varphi_w'(t\lambda)||\mathrm{d}(t\lambda)| \\\nonumber
   & \leq\frac{2}{3}\|g\|_{\infty}|\lambda|\int_0^1\frac{1-|t\lambda|^2}{(1-|t\lambda|)^2}\mathrm{d}t\\\nonumber
   &=\frac{2}{3}\|g\|_{\infty}|\lambda|\left(\frac{2}{|\lambda|}\log\frac{1}{1-|\lambda|}-1\right).
\end{align}
By combining (\ref{main-1}), (\ref{eq-march-18-1}) and (\ref {march-18-2}), we obtain
$$\left|(1-|z|^2)\Lambda _u(z)-(1-|w|^2)\Lambda _u(w)\right|\leq2\|g\|_{\infty}|\lambda|,$$
because
$$(1-|\lambda|^2)\left(\frac{2}{|\lambda|}\log\frac{1}{1-|\lambda|}-1\right)<2$$
for any $|\lambda|<1$.
This completes the proof. \qed
\section{Appendiex}\label{sec-app}
In this section, we calculate the precise value of $C_q$,
$$C_q=\left(\int_{\mathbb{C}}\frac{\mathrm{d}A(\xi)}{|\xi|^q|1-\xi|^q}\right)^{1/q},$$
where $1<q<2$. Note that in the proof of Theorem \ref{20-March-thm-1}, it was only required that this quantity is bounded.

Recall that the {\it hypergeometric function} $_pF_q$ is defined for $|z|<1$ by the power series (cf. \cite[(2.1.2)]{Landrews})
\be\label{bz-defn-2.1}_pF_q[a_1,a_2,\dots,a_p;b_1,b_2,\ldots,b_q;z]=\sum_{n=0}^\infty\frac{(a_1)_n\cdots(a_p)_n}{(b_1)_n\cdots(b_q)_n}\frac{z^n}{n!}.\ee
Here $(a)_n$ is the {\it Pochhammer symbol} and given as follows $(a)_n=\frac{\Gamma(n+a)}{\Gamma(a)}$.

\begin{Lem}\label{Lem-B}\textnormal{(cf. \cite[Theorem 2.1.2]{Landrews})}
The series $_{q+1}F_q(a_1,\ldots,a_{q+1};b_1,\ldots,b_q;x)$ converges absolutely for $|x|=1$ and \emph{{Re}}$(\sum_{m=1}^{q}b_m-\sum_{n=1}^{q+1}a_n)>0$. 
This series converges conditionally for $x=e^{i\theta}\neq1$ and $0\geq \emph{{Re}}(\sum_{m=1}^{q}b_m-\sum_{n=1}^{q+1}a_n)>-1$. This series diverges for \emph{{Re}}$(\sum_{m=1}^{q}b_m-\sum_{n=1}^{q+1}a_n)\leq -1$.
\end{Lem}

For $\eta=re^{it}\in\ID$, by using (\ref{bz-lem-2.0}) and Lemma \Ref{Lem-B}, we have
\begin{align*}
  \int_{\ID}\frac{\mathrm{d}A(\eta)}{|\eta|^q|1-\eta|^q} & =2\sum\limits_{n=0}^{\infty}\left(\frac{\Gamma(n+q/2)}{n!\Gamma(q/2)}\right)^2\int_0^1r^{2n+1-q}\mathrm{d}r \\
  &=\frac{_3F_2[1-\frac{q}{2}, \frac{q}{2}, \frac{q}{2};1, 2-\frac{q}{2};1]}{1-\frac{q}{2}}<\infty.
\end{align*}

\begin{figure}
\centering
\includegraphics[width=10cm]{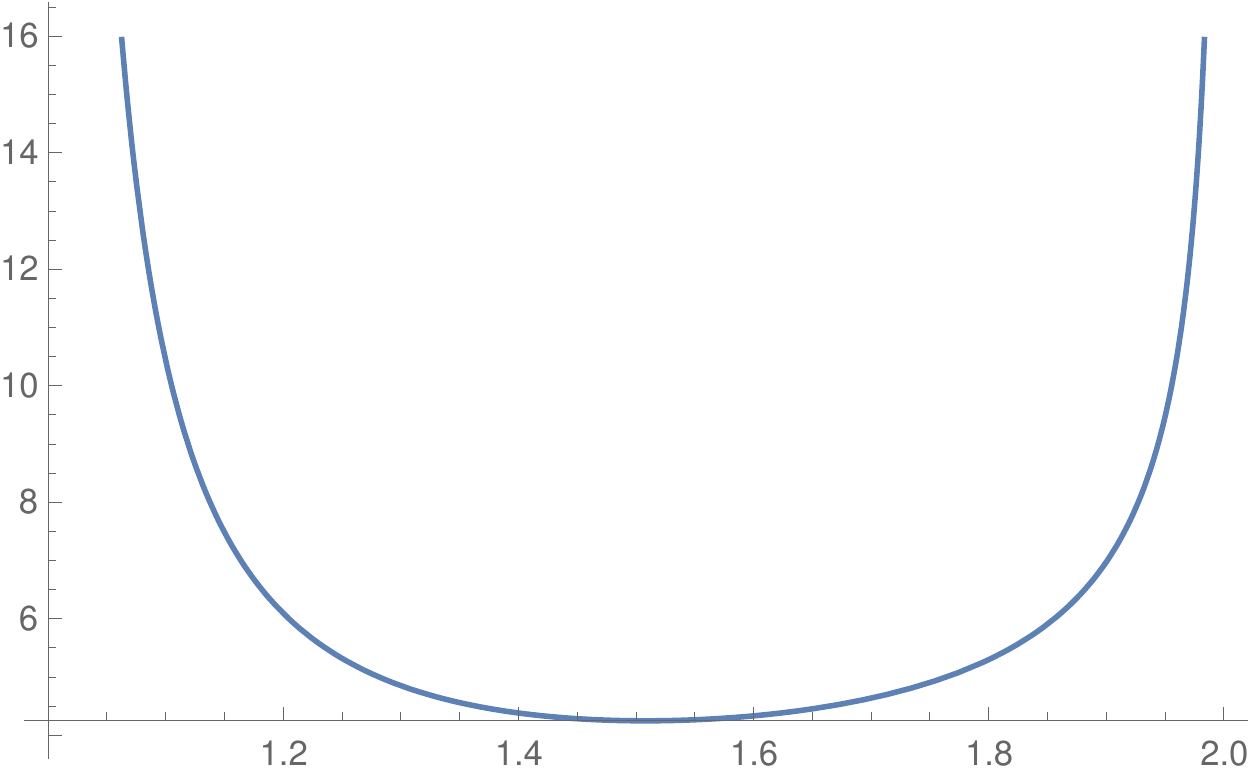}

\caption{The constant $C_q$ as a function of $q\in(1,2)$. }\label{Cqfig}
\end{figure}

For $\eta=re^{it}\in\mathbb{C}\setminus\ID$, again by (\ref{bz-lem-2.0}) and Lemma \Ref{Lem-B}, we have
\begin{align*}
\int_{\mathbb{C}\setminus\ID}\frac{\mathrm{d}A(\eta)}{|\eta|^q|1-\eta|^q} &=\int_{\ID}\frac{|\eta|^{2q-4}}{|1-\eta|^q}\mathrm{d}A(\eta) \\
  &=\frac{_3F_2[q-1, \frac{q}{2}, \frac{q}{2};1, q;1]}{q-1}<\infty.
\end{align*}
By combining the above two identities, we obtain the constant $C_q$. Numerical values of $C_q$ for $q\in(1,2)$ are illustrated in Figure \ref{Cqfig}.

\vspace*{5mm}
\noindent {\bf Acknowledgments}.
The research of the authors were supported by NSFs of China (No. 11501220, 11971124, 11971182), NSFs of Fujian Province (No. 2016J01020, 2019J0101), Subsidized Project for Postgraduates' Innovative Fund in Science Research of Huaqiao University
and the Promotion Program for Young and Middle-aged Teachers in Science and Technology Research of Huaqiao University (ZQN-PY402).

\end{document}